\documentclass{amsart}
\usepackage{amsfonts,amsmath,amsthm,amssymb, amscd}

\newtheorem*{theorem}{Theorem}

\newtheorem{corollary}{Corollary}
\newtheorem{lemma}{Lemma}
\newtheorem*{example}{Example}

\DeclareMathOperator{\rank}{rank}
\def\gp#1{\langle #1 \rangle}
\def\m1{^{-1}}

\author{V.~Bovdi, \fbox{V.~Rudko}}

\title[irreducible ring  representation]{On the  irreducible   representation algebra  of the alternating group of degree four}

\address{Institute of Mathematics, University  of Debrecen, H-4010 Debrecen, P.O.B. 12, Hungary}
\email{vbovdi@math.klte.hu}

\thanks{The research was supported by OTKA  No.K68383}

\subjclass {Primary: 20C11; Secondary:  16A64, 20C10,  20G05}

\keywords{integral representation, algebra representation}

\begin{document}

\begin{abstract}
We obtain a description  of the  irreducible  representation  algebra of the alternating group of degree four over the ring of $2$-adic integers.
\end{abstract}
\maketitle

Let $G$ be a finite group and let $K$ be a commutative   principal ideal domain. We denote the  module  $M$ of a   $K$-representation of the group  $G$ by  the symbol  $[M]$. We assume that $[M]=[N]$ if and only if the  $KG$-modules $M$ and  $N$ of $K$-representations are isomorphic. The ring  $\mathfrak{a}(KG)$ of $K$-representation of the group  $G$ is a smallest ring which contains all symbols $[M]$ with the following operations
\[
[M]+[N]=[M\oplus N],\qquad  [M]\cdot [N]=[M\otimes N],
\]
where\quad  $g(m\otimes n)=g(m)\otimes g(n)$, \quad  ($g\in G$,  $m\in M$, $n\in N$).

The subring  $\mathfrak{b}(KG)\subset \mathfrak{a}(KG)$ which is generated by  $[M]$,  where  $M$ is a module of the irreducible  $K$-representation of the group  $G$  is called the irreducible $K$-representations ring of  $G$.

The  algebras\quad  $\mathfrak{A}(KG)=\mathfrak{a}(KG)\otimes_{\mathbb{Z}}\mathbb{Q}$ and
$\mathfrak{B}(KG)=\mathfrak{b}(KG)\otimes_{\mathbb{Z}}\mathbb{Q}$
over the field of rational numbers $\mathbb{Q}$ are called the $K$-representation algebra of $G$ and the irreducible $K$-representation algebra of $G$, respectively.

It is well known that $\mathfrak{A}(RG)$, where $R$ is  the ring of  $p$-adic integers,  is finite dimensional if and only if  the $p$-Sylow subgroups of $G$ are cyclic of order $p^r$, where $r\leq 2$.   This   is a classical result of S.D.~Berman, P.M.~Gudivok, A.~Heller and I.~Reiner (see Theorem 33.6, \cite{Curtis_Reiner_book}, p.690).

Although the number of irreducible $R$-representations of the group $G$ is finite, the algebra $\mathfrak{B}(RG)$  may be infinite dimensional (see \cite{Barannik_Gudivok_Rudko, Gudivok_Gonchar_Rudko, Gudivok_Rudko_book}).

In the present paper we obtain  a description of the   irreducible  $R$-representation  algebra of the alternating group  $A_4$ of degree four over the ring $R$ of $2$-adic integers.  Note that earlier it was known that this algebra is infinite dimensional (see Proposition 13.1, \cite{Gudivok_Rudko_book}, p.108).

Our main result is the following.

\begin{theorem} Let  $R$ be the ring of  $2$-adic integers.
The   irreducible  $R$-representation algebra $\mathfrak{B}(RG)$ of the alternating group $G$ of degree $4$  isomorphic to
\[
\mathbb{Q}[x,x\m1]\oplus \mathbb{Q}[x,x\m1]\oplus \mathbb{Q}[x,x\m1]\oplus \mathbb{Q}[x,x\m1].
\]
\end{theorem}

\bigskip

\noindent
\underline{Notation}. Let $G$ be   the alternating group of degree four. Put  $a_1=(1,2)(3,4)$,\quad  $a_2=(1,4)(2,3)$,\quad  $b=(1,2,3)\in S_4$. Clearly
\[
G\cong \gp{a_1,  b}=H\rtimes \gp{b}\cong (C_2\times C_2)\rtimes C_3,
\]
 where
$H=\gp{a_1,a_2}\cong C_2\times C_2$ is the $2$-Sylow subgroup of $G$.

Let  $L$  be a free   $RG$-module over the ring $R$  with $R$-basis $\{e_1,e_2,e_3, e_4\}$
in which $G$ acts as a permutation group. Let  $d$ be  a  divisor of $4$ and put
\[
\begin{split}
L_d&=\Big\{\; \sum_{j=1}^4\alpha_je_j \;\mid\;  \alpha_j\in R,\quad  \sum_{j=1}^4\alpha_j\equiv 0\pmod{d}\; \Big\},\\
L_0&=R(e_1+\cdots+e_4), \qquad  M_d=L_d/L_0.\\
\end{split}
\]
Obviously $L_d$,  $L_0$  and  $M_d $  are $RG$-modules and the elements
$$
u_1=de_1+L_0, \qquad  u_2=e_2-e_1+L_0, \qquad  u_3=e_3-e_2+L_0
$$
form an $R$-basis of  $M_d$. Since  $e_4\equiv -e_1-e_2-e_3\pmod{L_0}$, it is not difficult to show that  the following
$R$-representation of  $G$
$$
\Gamma_d:a_1 \mapsto
\left(\begin{smallmatrix}
1&0&-4d^{-1}\\
d&-1&-2\\
0&0&-1
\end{smallmatrix}\right), \qquad
a_2\mapsto
\left(\begin{smallmatrix}
-3&4d{-1}&0\\
-2d&3&0\\
-d&2&-1
\end{smallmatrix}\right), \qquad
b\mapsto
\left(\begin{smallmatrix}
1&0&0\\
d&0&-1\\
0&1&-1
\end{smallmatrix}\right)
$$
is afforded by module the $M_d$.

Note that  $\Gamma_2$ is equivalent to the following monomial representation
$$
a_1\mapsto
\left(\begin{smallmatrix}
-1&0&0\\
0&1&0\\
0&0&-1
\end{smallmatrix}\right),\qquad\quad
b\mapsto
\left(\begin{smallmatrix}
0&0&1\\
1&0&0\\
0&1&0
\end{smallmatrix}\right).
$$

\bigskip

We shall need the following result.
\begin{lemma}\label{L:1}
(see \cite{Nazarova_A4})  The representations $\Gamma_d$, where $d=1,2,4$,
are irreducible and nonequivalent  $R$-represen\-ta\-tions of the group $G$. Moreover, except these representations and the trivial one
$\tau_0: g\mapsto 1$  ($g\in G$), the group $G$ has only one  more irreducible  $R$-representation:
$$
\tau:a_1\mapsto  E,\qquad  b\mapsto \left(\begin{smallmatrix}
0&-1\\
1&-1
\end{smallmatrix}\right).
$$
\end{lemma}
\begin{lemma}\label{L:2}
Let $\Gamma=P_0\oplus P_1$, where $P_0=RGw_0$,  $P_1=RGw_1$ are projective  $RG$-modules and
$w_0=\frac{1}{3}(1+b+b^2)$,  $w_1=1-w_0$ are orthogonal idempotents. Then the  following equations hold
\[
[P_0]^2=2[P_0]+[P_1], \qquad  [P_0][P_1]=2[P_0]+3[P_1], \qquad [P_1]^2=6[P_0]+5[P_1].
\]
\end{lemma}

\begin{proof}
It is easy to see that  $\chi_{P_0}(b)=1$ and $\chi_{P_1}(b)=-1$, where $\chi_{P_{i}}$ is the character  of $P_{i}$. Remember that
$\chi_{P_0\oplus P_1}=\chi_{P_0}+\chi_{P_1}$, \quad  $\chi_{P_0\otimes P_1}=\chi_{P_0}\chi_{P_1}$.

If $[P_0]^2=s[P_0]+t[P_1]$, then
$s-t=1$ and   $4s+8t=16$. It follows  that $s=2$ and  $t=1$, so the first equation is true. The proof of second one is analogues.
\end{proof}

\begin{corollary} \label{C:1}
The elements
\begin{equation}\label{E:1}
\mathfrak{f}_1=\textstyle\frac{1}{12}[\Gamma],\qquad \mathfrak{f}_2=[P_0](1-\mathfrak{f}_1),\qquad \mathfrak{f}_3=(1-[P_0])(1-\mathfrak{f}_1)
\end{equation}
are pairwise orthogonal idempotents of the algebra  $\mathfrak{A}(RG)$  such that
$\mathfrak{f}_1+\mathfrak{f}_2+\mathfrak{f}_3=1$. Moreover, \qquad
\[
\mathfrak{A}(KG)\mathfrak{f}_1 \cong \mathfrak{A}(KG)\mathfrak{f}_2 \cong\mathbb{Q}.
\]
\end{corollary}

Let  $W$ be a finite group.  Assume that for the $K$-representations of  $W$ the Krull-Schmidt-Azumaya theorem holds (see Theorem 36.0, \cite{Curtis_Reiner_book}, p.768). We assume that there exists  at least one nonzero indecomposable, non-projective $KW$-module of the  $K$-representation of the group $W$. For each module $B$ of  the  $K$-representation of the group $W$ we consider
the  $K$-module $B^*=Hom_K(B,K)$. Clearly $B^*$ is a $KW$-module of  the  $K$-representation of the group $W$ such that
\[
(g\cdot \varphi)\cdot b=\varphi(g\m1b), \qquad\quad  (g\in W, \;  \varphi\in B^*,\;  b\in B).
\]
The module $B^*$ is called the contragredient module of $B$.

Note that each module of  the $K$-representation of $W$ is always a homomorphic image  of a free module,
and  it can be considered as a submodule of a free module, too.

We shall use  the following well known result.
\begin{lemma}\label{L:3}
Let\quad $0\rightarrow B\rightarrow P\rightarrow A\rightarrow 0$\quad  be an exact sequence of modules of $K$-representations of the  group $W$.  If  $P$ is a projective  $KW$-module and   $A$ is an indecomposable non-projective $KW$-module,
then the module  $B$ contains only one unique indecomposable  nonprojective  direct summand  $B_0$.  Moreover the  sequence
$
0\rightarrow B_0 \rightarrow P_0 \rightarrow A \rightarrow 0
$
\quad
is exact for some projective  $KW$-module $P_0$.
There  exists a duality corresponding to the exchange of  $A$ and $B$.
\end{lemma}
\begin{proof} All modules in the lemma are free of finite rank over the ring  $K$.
There exists a projective  $KW$-module  $P$ and a sequence of  indecomposable  $KW$-modules $B_1, \ldots, B_s$,
such that the sequence\quad
$0\rightarrow B_1\oplus \cdots \oplus B_s \rightarrow P
\rightarrow A \rightarrow 0$\quad
is exact.

Assume that   $A_1, \ldots, A_s$ are  $KW$-modules  and  $P_1, \ldots, P_s$ are projective  $KW$-modules,
such that the  sequences\quad
$0\rightarrow B_j \rightarrow P_j \rightarrow A_j \rightarrow 0$\quad are exact for each $1\leq j\leq s$.
Sum them over  $j$ and apply Schanuel's Lemma (Lemma 2.24, \cite{Curtis_Reiner_book}, p.30), then    we obtain that
\[
P\oplus A_1\oplus \cdots  \oplus A_s \cong A\oplus P_1 \oplus  \cdots \oplus
P_s.
\]
This yields  that $A$ is a direct summand  only  one of the $A_i$, says of $A_1$, so  $A_2,  \ldots ,A_s$ are projective.
Then the modules  $B_2, \ldots ,B_s$ are projective,  too. Therefore we have the following exact sequence
\quad $0\rightarrow B_0\oplus P' \rightarrow P \rightarrow A
\rightarrow 0$,\quad
where  $P'$ is projective and  $B_0$ is an indecomposable non-projective  $KW$-module.
The following isomorphism of  $KW$-modules hold
$$
A\cong P/(B_0 \oplus P')\cong (P/P')/B_0,
$$
where  $P/P'$ is a projective $KW$-module.

By the use of the  contragredient module we obtain the dual statement. \end{proof}

\bigskip

Obviously,   each    projective  $RH$-module is free, where $H=Syl_2(G)$.  Let  $M$ be a module of
$R$-representation of the group  $H$ and  assume that $M$ does not contain  projective summands.
By  $\Theta(M)$ we denote the kernel (the Green operator, after J.A.~Green) of the projective cover   of $RH$-module $M$
(see Lemma \ref{L:3}). The following sequence of $RH$-modules
$$
0\rightarrow \Theta(M)\rightarrow F\rightarrow M\rightarrow 0
$$
is exact, where    $F$ is a free $RH$-module of smallest rank. If $M$ is an indecomposable  $RH$-module, then
$\Theta(M)$ is indecomposable, too.

Let  $\Delta_0$ be the module of the trivial representation\;   $h\mapsto 1$ \;  ($h\in H$) of $H$. Put
\[
\Delta_n=\Theta(\Delta_{n-1}), \qquad  \Delta_{-n}=\Delta_n^*,  \qquad
(n=1,2, \ldots)
\]
where $\Delta_n^*$ is the contragredient  module to  $\Delta_n$.

\begin{lemma} \label{L:4}
The  sequence\quad $
0\rightarrow \Delta_n \rightarrow (RH)^n \rightarrow
\Delta_{n-1}\rightarrow 0$\quad
is exact and
$$
\rank_R(\Delta_n)=2n+1.
$$
\end{lemma}

\begin{proof}  Let $\Gamma^n$ be a $n$-dimensional module  over the ring $RH$, where $H=\gp{a_1,a_2}$. In the follows we will treat
$\Gamma^n$  considering as component columns.
Denote by  $F_n$  the  matrix over the ring  $RH$ with  $n$ rows and   $n+1$ columns, in which only the
three upper diagonals (if we begin calculation from the main diagonal) are non-zero:
\begin{itemize}
\item[$\bullet$] the main diagonal of $F_n$ is the first $n$ element from the following sequence
$\digamma_1\cup \digamma_1\cup \digamma_1,\ldots$, where
\[
\digamma_1=
\begin{cases}
\{a_1-1, a_2-1, -(a_1+1), -(a_2+1)\} & \quad \text{if} \quad n\quad \text{is odd},\\
\{a_1+1, a_2+1, -(a_1-1), -(a_2-1)\} & \quad \text{if} \quad n\quad \text{is even}.\\
\end{cases}
\]
\item[$\bullet$]  the second diagonal consist  of  zero elements except the last one which is equal to one of the following elements
\[
a_2-1,\;  -(a_1-1),\;   a_2+1,\;  a_1+1
\]
according to the cases  $n\equiv \{1,2,3,0\}\pmod{4}$, respectively.
\item[$\bullet$]  the third  diagonal consist of  $n-1$ elements of the sequence
\newline
$\digamma_2\cup \digamma_2\cup \digamma_2,\ldots$, where  $\digamma_2=\{a_2-1, a_1-1, a_2+1, a_1+1\}$.
\end{itemize}

\begin{example}
\qquad\quad
$F_4=
\left(
\begin{smallmatrix}
a_1+1 & 0     & a_2-1   &    0       &      0    \\
 0    & a_2+1 & 0       &  a_1-1     &      0    \\
 0    & 0     & -(a_1-1)&    0       &   a_2+1   \\
 0    & 0     & 0       &   -(a_2-1) &   a_1+1   \\
\end{smallmatrix}
\right);
$
\[
F_5=
\left(
\begin{smallmatrix}
a_1-1 & 0     & a_2-1   &    0       &      0 &   0    \\
 0    & a_2-1 & 0       &  a_1-1     &      0 &   0    \\
 0    & 0     & -(a_1+1)&    0       &  a_2+1 &   0    \\
 0    & 0     & 0       &   -(a_2+1) &    0   &  a_1+1 \\
 0    & 0     & 0       &   0        &  a_1-1 &  a_2-1 \\
\end{smallmatrix}
\right).
\]
\end{example}
Let $V_n$ be a submodule in  $\Gamma^n $ generated  by the columns of the matrix $F_n$.  Let $\varepsilon: \Gamma^n
\rightarrow V_{n-1}$ be the following epimorphism  of modules:
$$
\Gamma^n \ni x=(x_1, \ldots, x_n)\mapsto  x_1f_1+ \cdots +x_nf_n \in V_n,
$$
where $f_j$  are columns of the matrix $F_{n-1}$.

It is not difficult  to check that \qquad $F_{n-1} F_n=0$.

Moreover, any solution of $\varepsilon (x)=0$ belongs to the linear combination of the columns of
$F_n$, so  $x \in V_n$. Therefore, we proved the exactness  of  the sequence
$$
0\rightarrow V_n \rightarrow \Gamma^n \rightarrow V_{n-1}
\rightarrow 0.
$$
We construct a basis of the $R$-module $V$ in the following way: the first  $n+1$ elements of the
basis are the columns   of  $F_n$. Then we take the first $n$ columns of $F_n$, which multiply
each of them by corresponding  element of  the sequence
$ \digamma_3\cup \digamma_3\cup\ldots \cup\digamma_3$, where
$\digamma_3=\{a_2-1, a_1-1, a_2+1, a_1+1\}$. These columns form the remaining $n$ elements of the basis.

Consequently  we obtain a system of  $2n+1$ basis elements of  $R$-module $V_n$.
By Lemma \ref{L:3}, beginning with  $V_0$  all  $RH$-modules  $V_n$ are indecomposable. \end{proof}

Induce the  exact sequence of  Lemma \ref{L:4}  from the subgroup $H$ to the group  $G$. In this exact sequence
of  $RG$-modules  the middle  modules are free, and  the $2^{nd}$  and the $4^{th}$ terms of the sequence  have a decomposition
into the direct sums of indecomposable  $RG$-modules
$$
\Delta_n^G=\Delta_{0,n}\oplus \Delta_{1,n},
$$
where\; $\rank(\Delta_{0,n})=\rank(\Delta_n)$\; and \;  $\rank(\Delta_{1,n})=2\cdot \rank(\Delta_n)$.
Moreover  $\Delta_{0,0}$ and $\Delta_{1,0}$ are modules of the  $R$-representations  $\tau_0$, $\tau$ (see Lemma \ref{L:1}) of $G$, respectively.
It is easy to check that $\Delta_{0,1}$ and $\Delta_{0,-1}$ are modules of the irreducible $R$-representations of  $\Gamma_1$,
$\Gamma_4$ of  $G$, respectively (see Lemma  \ref{L:1}). Since
\[
(\Delta_n^G)|_H=(\Delta_n)^{(3)},
\]
the $RG$-module  $\Delta_{0,n}$ is a lifting  of the $RH$-module $\Delta_n$.

We will use  the notation $\Delta_n=\Delta_{0,n}$.

\begin{lemma} \label{L:5}
The sequences of  $RG$-modules
\[
\begin{split}
0&\rightarrow  \Delta_{3k}\rightarrow  \Gamma^{k} \rightarrow \Delta_{3k-1}\rightarrow  0,\\
0&\rightarrow  \Delta_{3k+1}\rightarrow  P_0\oplus\Gamma^{k} \rightarrow  \Delta_{3k}\rightarrow  0,\\
0&\rightarrow  \Delta_{3k+2}\rightarrow  P_1\oplus\Gamma^{k} \rightarrow  \Delta_{3k+1}\rightarrow  0\\
\end{split}
\]
are exact, where $\Gamma=RG=P_0\oplus P_1$ (see Lemma \ref{L:2}). Moreover
\[
\begin{split}
\Delta_{3k}\otimes \Delta_1   &=\Delta_{3k+1}\oplus \Gamma^{k},\\
\Delta_{3k+1}\otimes \Delta_1 &=\Delta_{3k+2}\oplus P_0 \oplus\Gamma^{k},\\
\Delta_{3k+2}\otimes \Delta_1 &=\Delta_{3k+3}\oplus  P_1 \oplus\Gamma^{k}.\\
\end{split}
\]
\end{lemma}

\begin{proof}  It is easy to see that
\[
\{ (a_1-1)(a_2-1)w_0, \, (a_1-1)w_0, \, (a_2-1)w_0\}
 \]
is an  $R$-basis of the  $RG$-submodule  $M\subset P_0$. Obviously, $P_0 /M\cong \Delta_0$,
so we can assume that  $M=\Delta_1$ and  the following sequence
\begin{equation}\label{E:2}
0\rightarrow  \Delta_1\rightarrow  P_0 \rightarrow   \Delta_0
\rightarrow  0
\end{equation}
is exact.
Comparing the values of the characters at  $b\in G$, we get
\[
P_0\otimes \Delta_1 \cong \Gamma=P_0\oplus P_1.
\]
Multiplying (\ref{E:2})  tensorially  by  $\Delta_1$, we obtain the exact sequence
$$
0\rightarrow  \Delta_1 \otimes \Delta_1\rightarrow  P_0\oplus P_1 \rightarrow
\Delta_1 \rightarrow  0,
$$
 which is  possible only if $\Delta_1\otimes\Delta_1 \cong \Delta_2\oplus P_0$, so the sequence
$$
0\rightarrow  \Delta_2\rightarrow  P_1
\rightarrow   \Delta_1 \rightarrow  0.
$$
is exact, too.
Multiplying  the last sequence again  tensorially by  $\Delta_1$ and using the value of the characters
(for example, $\chi_{\Delta_2}(b)=-1$),  we can verify   the lemma  for   $n=2,3,\ldots$.
\end{proof}

Using  the contragredient modules we can obtain an  analogue of Lemma \ref{L:5}  for the modules $\Delta_m$ with negative $m$.
As a consequence we have

\begin{corollary} \label{C:2}
In the algebra  $\mathfrak{A}(RG)$ the following equations hold
$$
[\Delta_n]\cdot [\Delta_m]\mathfrak{f}_3=[\Delta_{n+m}]\mathfrak{f}_3, \qquad [\Delta_{1,0}]\cdot [\Delta_{1,0}]=2[\Delta_0]+[\Delta_{1,0}].
$$
\end{corollary}

\begin{proof} These equations follows from Lemma \ref{L:5}, where  $[P_j]\mathfrak{f}_3=0$ and $\mathfrak{f}_3$ from (\ref{E:1}).   \end{proof}

The group  $H=\gp{a_1,a_2}\cong C_2\times C_2$ has the following four linear characters:
\[
\begin{matrix}
\delta_0: &a_1\mapsto  1, & a_2 \mapsto & \;\,1;  & \quad &\delta_1:a_1\mapsto  -1, & a_2 \mapsto  \;\, 1;\\
\delta_2: &a_1\mapsto  1, & a_2 \mapsto &-1;  &\quad  &\delta_3:a_1\mapsto  -1, & a_2 \mapsto  -1.
\end{matrix}
\]
It is easy to check  that the induced representations $\delta_1^G$,  $\delta_2^G$,  $\delta_3^G$  of  $G$ are irreducible  and equivalent to the representation  $\Gamma_2$ (see Lemma \ref{L:1}).

\begin{lemma}\label{L:6}
The following equations hold
\[
[\Delta_{1,0}]^2=[\Delta_{1,0}]+[\Delta_{0}],\qquad [L]^2=[\Delta_0]+[\Delta_{1,0}]+2[L],\\
\]
where $L$ is the module of the representation $\Gamma_2$ (see Lemma \ref{L:1}) and $\Delta_{0}$ is the trivial $RG$-module.
\end{lemma}

\begin{proof}
Using  Mackey's theorem (see Theorem 10.18, \cite{Curtis_Reiner_book}, p.240) we get
\[
\delta_2^G\otimes \delta_2^G\simeq \sum_{j=0}^{2}(\delta_2\otimes
\delta_2^{b^j})^G,
\]
where
$\delta_2^{b^j}(h)=\delta_2(b^{-j}hb^j)$ for   $h\in H$. This  yields   that
$
\delta_2\otimes \delta_2=\delta_0$,\quad  $\delta_2\otimes \delta_2^{b}=\delta_1$,\quad
$\delta_2\otimes \delta_2^{b^2}=\delta_3$,\quad  and \quad $\delta_2^G\otimes \delta_2^G=   \delta_0^G +\delta_1^G+\delta_3^G$. \end{proof}

\begin{proof}[Proof of the Theorem] Let  $\mathfrak{B}'(RG)=\mathfrak{B}(RG)\mathfrak{f}_3$.  Put
\[
x=[\Delta_1],\qquad    y=[\Delta_{1,0}], \qquad  z=[L].
\]
We identify  $\mathfrak{f}_3$ with the unity  $1$ of the algebra $\mathfrak{B}'(RG)$, where $\mathfrak{f}_3$ from (\ref{E:1}). Since $[\Delta_1][\Delta_{-1}]\mathfrak{f}_3=\mathfrak{f}_3$, we can assume that $[\Delta_{-1}]\mathfrak{f}_3=\frac{1}{x}$.  By Lemma \ref{L:6}  it follows  that
\[
\mathfrak{B}'(RG)=\gp{\; 1, \; x,\; \textstyle\frac{1}{x}, \; y,\; z \; \mid
\;  y^2=y+2, \quad  z^2= 2z+y+1\;}.
\]
Since $y^2-y-2=(y-2)(y+1)$, we have
\[
\gp{\; 1, \; x,\; \textstyle\frac{1}{x}, \; y\;}\cong \mathbb{Q}[x,\frac{1}{x}][y]/\gp{y^2-y-2}\cong \mathbb{Q}[x,\frac{1}{x}]\oplus \mathbb{Q}[x,\frac{1}{x}].
\]
Now  the algebra $\mathfrak{B}'(RG)\cong
\mathbb{Q}[x,\frac{1}{x}]\oplus \mathbb{Q}[x,\frac{1}{x}]\oplus \mathbb{Q}[x,\frac{1}{x}]\oplus \mathbb{Q}[x,\frac{1}{x}]$, because
\[
z^2-2z-y-1=
\begin{cases}
(z+1)(z-3)\qquad &\text{for} \qquad y=2;\\
(z-2)z \qquad &\text{for} \qquad y=-1.\\
\end{cases}
\]
Finally, since the algebra  $\mathfrak{B}(RG)$ has no projective summands, the map $u\mapsto u\mathfrak{f}_3$ gives the isomorphism
$\mathfrak{B}'(RG)\cong \mathfrak{B}(RG)$.
\end{proof}

\newpage

\bibliographystyle{abbrv}
\bibliography{Algebra_A4_English}

\end{document}